\documentclass[reqno]{amsart}

\usepackage{amsfonts}
\usepackage{amssymb}
\usepackage{amsmath}
\usepackage[all]{xy}

\newenvironment{definition-2}{\textbf{Definition 2.}}{}
\newenvironment{example-3}{\textbf{Example 3.}}{}
\newenvironment{example-4}{\textbf{Example 4.}}{}
\newenvironment{theorem-8}{\textbf{Theorem 8.}}{}
\newenvironment{theorem-9}{\textbf{Theorem 9.}}{}
\newenvironment{theorem-11}{\textbf{Theorem 11.}\it}{\normalfont}

\renewcommand{\mapsto}{\longmapsto}
\renewcommand{\to}{\longrightarrow}

\newcommand{\A}{\ensuremath{\underline{A}}}
\newcommand{\B}{\ensuremath{\underline{B}}}
\newcommand{\X}{\ensuremath{\underline{X}}}
\newcommand{\Y}{\ensuremath{\underline{Y}}}
\renewcommand{\S}{\ensuremath{\underline{S}}}
\newcommand{\SSigma}{\ensuremath{\underline{\Sigma}}}
\newcommand{\Sets}{\ensuremath{\underline{\mathrm{Sets}}}}
\newcommand{\cls}{\ensuremath{\mathrm{cls}}}

\begin{document}

\newdir{>}{{}*:(1,-.2)@^{>}*:(1,+.2)@_{>}}
\newdir{<}{{}*:(1,+.2)@^{<}*:(1,-.2)@_{<}}

\title{On satellites in arbitrary categories}
\author{G.\ Janelidze}
\maketitle

\vspace{-3ex}
\begin{center}
(Communicated by Academician G.\ Chogoshvili\\
on 20 December 1975)\footnote{Originally published in Russian as G.~Z.~Janelidze, \emph{On satellites in arbitrary categories}, Bull.\ Georgian Acad.\ Sci.\ \textbf{82} (1976), no.~3, 529--532. Translated by Jone Intxaurraga Larra\~naga and Tim Van der Linden with the help of Alexander Frolkin and Julia Goedecke.}
\end{center}

\bigskip
We generalize the definition of satellites with respect to presheaves (and copresheaves) with trace in the sense of \cite{Inasaridze75}; a presheaf with trace is replaced by a graph with a pair of diagrams defined on it.

We show that the right satellite functor is left adjoint to the
left satellite functor, and that a functor having a right (left)
adjoint preserves right (left) satellites.

In particular cases the construction of satellites is given.

We shall denote by $\S(\X,\Y)$, where $\X$ and $\Y$ are
categories, the category of which the objects are graphs $\S$
together with two diagrams $F=F(\S)\colon \S\to \X$ and
$G=G(\S)\colon \S \to \Y$; the morphisms in $\S(\X,\Y)$ are
defined in the natural way.

We have an isomorphism
\begin{equation}
\S(\X,\Y)^{\circ}\approx\S(\Y^{\circ},\X^{\circ})
\end{equation}
given by $(\xymatrix@1{\X & \S \ar[l]_-{F} \ar[r]^-{G} &
\Y})\mapsto (\xymatrix@1{\Y^{\circ} & \S^{\circ} \ar[l]_-{G}
\ar[r]^-{F} & \X^{\circ}})$.

\begin{definition-2}
In the diagram
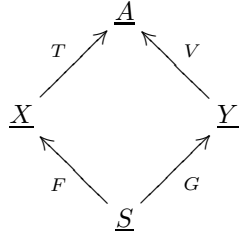
\[
\xymatrix@!{ & {\A} \\
\X \ar[ru]^-{T} && \Y \ar[lu]_-{V}\\
& \S \ar[lu]^-{F} \ar[ru]_-{G}}
\]
of graphs and their morphisms, let the graphs $\X$, $\Y$ and $\A$
be categories, and let the diagrams $T$ and $V$ be functors. We
will call a morphism of diagrams $\delta\colon{TF\to VG}$ an
\emph{$\S$-connecting morphism} (more precisely, an
\emph{$(\S,F,G)$-connecting morphism}) $\delta\colon{T\to V}$. In
this case we shall call the triple $(T,\delta,V)$ an
\emph{$\S$-connected pair (of functors with values in $\A$).} If
an $\S$-connected pair $(T,\delta,V)$ is right universal (left
universal), then we say that $V$ is a \emph{right satellite} for
$T$ ($T$ is a \emph{left satellite} for $V$) relative to $\S$ and
we write $V=\S^{1}T$ ($T=\S_{1}V$).
\end{definition-2}

The satellites of functors with respect to copresheaves (with
trace) in the sense of~\cite{Inasaridze75} are special cases of
satellites in the sense of Definition~2. This also applies to
satellites of contravariant functors relative to presheaves.
Moreover, unlike~\cite{Inasaridze75}, Definition 2 allows to
consider satellites of (covariant) functors relative to presheaves
and satellites of contravariant functors relative to copresheaves.

\begin{example-3}
In Situation~2, let $\X=\Y$, let $\S$ be a subgraph of the
category of diagrams $\X^{\SSigma}$, where $\SSigma$ is the graph
of the form $\SSigma = (\xymatrix@1{{\Sigma''} \ar[r] &{\Sigma}
\ar[r] & {\Sigma'}})$, and let $F(\S)$ and $G(\S)$ be defined by
$F(\S)(-)=(-)(\Sigma')$ and $G(\S)(-)=(-)(\Sigma'')$. If now $\S $
is a full subcategory of $\X^{\SSigma}$ and $|\S|$ satisfies the
conditions of \cite[Def.\ 1.3]{Inasaridze75}, then Definition~2
coincides with \cite[Def.\ 1.5, 1.6]{Inasaridze75} and, in
particular, is a generalization of the classical definition.
\end{example-3}

\begin{example-4}
Let $K\colon \Y\to \A$ be a functor and $\SSigma$ a graph with a
distinguished object $\Sigma$. We shall denote by $\SSigma'$ the
full subgraph of $\SSigma$ such that
$|\SSigma'|=|\SSigma|-\{\Sigma\}$. We will furthermore choose: a
subgraph $\S$ of the category $\Y^{\SSigma}$, whose objects will
be called \emph{resolutions}; a subcategory $\X$ of the category
$\A^{\SSigma'}$, whose objects will be called \emph{complexes};
and a functor $T\colon \X\to \A$, which will be called the
\emph{homology} (or \emph{homotopy}) functor. We shall require the
functor $K$ to induce a diagram $\S\to \X$, which we shall denote
by $F(\S)$. Next we define the diagram $G(\S)\colon \S\to \Y$ by
$G(\S)(-)=(-)(\Sigma)$. The functor $\S^{1}T\colon \Y\to \A$
(provided it exists) will be called \emph{the derivative of the
functor $K$, relative to the triple $(\S,\X,T)$}. Under an
appropriate choice of this triple, one obtain the derivatives in
the sense of \cite{Inasaridze70} (from which the idea of
considering derived functors as homologies of satellite functors
is also taken) and \cite{Dold-Puppe, Eilenberg-Moore,
Tierney-Vogel2, Swan}; here it makes no difference whether left or
right derivatives are used.
\end{example-4}

In Examples 3 and 4, only the covariant case was considered; the
contravariant case can be considered analogously.

We return again to Situation 2 and will suppose for simplicity
that $\A$ admits satellites relative to $(\X,\Y)$, i.e.\ $\S^{1}T$
and $\S_{1}V$ exist for all $\S\in\S(\X,\Y)$, $T\colon{\X\to \A}$
and $V\colon {\Y\to \A}$. Then the assignments $(\S,T)\mapsto
\S^{1}T$ and $(\S,V)\mapsto\S_{1}V$ determine bifunctors
\setcounter{equation}{4}
\begin{gather}
\S(\X,\Y)\times(\X,\A)\to(\Y,\A),\\
\S(\X,\Y)^{\circ}\times(\Y,\A)\to(\X,\A),
\end{gather}
(where $(\X,\A)$ is the category of functors from $\X$ to $\A$ and
$(\Y,\A)$ has the analogous meaning) and there are natural
isomorphisms
\begin{equation}
\S^{1}T\approx (\S^{\circ})_{1}T,\quad \S_{1}V\approx (\S^{\circ})^{1}V
\end{equation}
(if we make $|(\X,\A)|=|(\X^{\circ},\A^{\circ})|$ and
$|(\Y,\A)|=|(\Y^{\circ},\A^{\circ})|$). In this way, we obtain the
bifunctors (5) and (6) from each other using (1).

Let us add two more properties of satellites:

\begin{theorem-8}
For each $\S\in \S(\X,\Y)$, the functor $\S^{1}\colon{(\X,\A)\to(\Y,\A)}$ is left adjoint to the functor $\S_{1}\colon{(\Y,\A)\to(\X,\A)}$.
\end{theorem-8}

\begin{theorem-9}
Let $K\colon{\A\to \B}$ be a functor. If the pair $(T,\delta,V)$ is right universal (left universal) and the functor $K$ has a right (left) adjoint, then the pair $(KT,K\delta,KV)$ is also right (left) universal.
\end{theorem-9}

We turn to the construction of satellites. Let $\A=\Sets$, $\X$
and $\Y$ small categories, and $T\colon{\X\to\A}$,
$\S\in\S(\X,\Y)$ fixed. We shall define a functor $\Y\cdot \S\cdot
T\colon{\Y\to\A}$ as follows. We put \setcounter{equation}{9}
\begin{equation}
\Y\cdot \S\cdot T(Y)=\frac{\displaystyle\coprod_{X\in|\X|,Y'\in|\Y|}\Y(Y',Y)\times (G^{-1}(Y')\cap F^{-1}(X))\times T(X)}{\sim},
\end{equation}
where $\sim$ is the smallest equivalence relation under which
\[
\cls(yG(s),S_{1},t)=\cls(y,S_{2},TF(s)(t))
\]
for each $\S$-morphism $s\colon{S_{1}\to S_{2}}$. If $f\in \Y$,
then we define $\Y\cdot \S\cdot T(f)$ by $\cls(y,S,t)\mapsto
\cls(fy,S,t)$.

\begin{theorem-11}
Let $\A=\Sets$ and let $\X$, $\Y$ be small categories, $\S\in\S(\X,\Y)$ and $T\colon{\X\to \A}$, $V\colon {\Y\to\A}$ functors. Then
\renewcommand{\theenumi}{\alph{enumi}}
\begin{enumerate}
\item the assignment
\[
\delta(S)\colon{t\mapsto \cls(1_{G(S)},S,t)}
\]
defines an $\S$-connecting morphism $\delta\colon{T\to \Y\cdot\S\cdot T}$ and the $\S$-connected pair $(T,\delta,\Y\cdot\S\cdot T)$ is right universal;
\item the assignment
\[
\vartheta(S)\colon{\varphi\mapsto \varphi(G(S))(\cls(1_{G(S)},S,1_{F(S)})}
\]
defines an $\S$-connecting morphism
\[
\vartheta\colon{(\Y,\A)(\Y\cdot\S\cdot(\X(-,?)),V(?))\to V(-)}
\]
and the $\S$-connected pair $((\Y,\A)(\Y\cdot\S\cdot(\X(-,?)),V(?)),\vartheta,V)$ is left universal.
\end{enumerate}
\end{theorem-11}

This also makes it possible to construct satellites, and in the
case when $\A$ is a variety of universal algebras (for instance,
the category of abelian groups), the passage from $\Sets$ to $\A$
works in the same way as it does for limits: for right satellites
one should take the $\A$-free algebra over $\Y\cdot\S\cdot T(Y)$
(for every $Y\in|\Y|$), then take the (minimal) quotient that
makes every $\delta(S)$ a homomorphism; for left satellites one
only needs to equip $(\Y,\A)(\Y\cdot\S\cdot(\X(X,?)),V(?))$ (for
every $X\in|\X|$) with an appropriate $\A$-algebra structure.

\bigskip
\begin{center}
A.\ Razmadze Mathematical Institute\\
Georgian Academy of Sciences\\
Tbilisi, Georgia\\
\smallskip
(Received 16 January 1976)
\end{center}

\newpage
\section*{Author's remarks, September 2008}

\small\noindent 1. Many thanks from the author to Tim Van der Linden, Julia Goedecke, and everyone who helped them in this translation!

\medskip
\noindent 2. The paper was written in Russian, because publishing a paper in English was not allowed in Soviet mathematical journals. Every Soviet author had ``two initials'' before his surname, the first of which was indeed initial while the second was his father's initial (``patronymic''), according to the Russian tradition, although it was inappropriate for Georgians (just as it would be inappropriate for Western-Europeans).

\medskip
\noindent 3. The author did not have his PhD yet, and so submitting a paper he needed an ``approval'' from a professor---in this case Hvedri Inassaridze, who was not yet an Academician (Soviet ``Academician'' means ``Member of Academy'', which is a title far above Professor, which itself is far above the PhD level)---and therefore his approval should have been followed by a ``final approval''---in this case by Academician George Chogoshvili.

\medskip
\noindent 4. There is a misprint in the formula (isomorphism) (\textbf{1}): the first ``$^{\circ}$'' should be removed (this obvious misprint was originally noticed by then young Beso Pachuashvili). In fact there is another duality there, but it is not really explored in the paper.

\medskip
\noindent 5. The original Russian version uses the term ``diagram scheme'' instead of ``graph''.

\medskip
\noindent 6. Most importantly, back in 1975 the author did not know about Kan extensions, but ``rediscovered'' them and their basic properties. In particular Theorem~11 follows from various known observations. A few years later the author saw ``All concepts are Kan extensions'' in Saunders Mac Lane's ``Categories for the working Mathematician'' (written in 1971); however, in 1990, Mac Lane suggested to him to publish his results on satellites fully. This has not happened so far, but some of those results are in

\smallskip
G.\ Janelidze, \emph{Satellites and Galois extensions of commutative rings}, PhD Thesis, Tbilisi, 1978

\smallskip
G.\ Janelidze, \emph{Satellites with respect to Galois extensions}, Proc.\ Math.\ Inst.\ Georgian Acad.\ Sci. \textbf{LXII} (1979), 38--48 (in Russian)

\smallskip
G.\ Janelidze, \emph{Computation of Kan extensions by means of injective objects and the functors Ext in nonadditive categories}, Proc.\ Math.\ Inst.\ Georgian Acad.\ Sci. \textbf{LXX} (1982), 42--51 (in Russian)

\smallskip
G.\ Janelidze, \emph{Internal categories and Galois theory of commutative rings}, DSc Thesis, Tbilisi, 1991 (defended at St.-Petersburg State University in 1992)

\smallskip
H.\ Inassaridze, \emph{Nonabelian homological algebra and its applications}, Mathematics and its Applications \textbf{421}, Kluwer Academic Publishers, Dordrecht, 1997.

\medskip
\noindent 7. The first sentence of the paper is not formulated very well: one can indeed get Inassaridze satellites as a special case, but what the author is doing is a ``better'' approach rather than a generalization. In fact the reason of ``better'' is that it is a situation where Grothendieck fibrations are better than indexed categories (which Inassaridze called ``presheaves of categories''), and moreover, fibrations can be replaced with arbitrary functors. And furthermore, things are ``inspired by Yoneda rather than Grothendieck'', but explaining all these is a long story... Note also that (nearly the same) description of Inassaridze satellites in terms of Kan extensions was obtained by Polish mathematician Stanislaw Balcerzyk in

\smallskip
S.\ Balcerzyk, \emph{On Inassaridze satellites relatively to traces of presheaves of categories}, Bull.\ Acad.\ Pol.\ Sci., S\'er.\ Sci.\ Math.\ Astron.\ Phys. \textbf{25} (1977), 857--861.


\begin{thebibliography}{1}

\bibitem{Inasaridze75}
H.\ Inassaridze, \emph{Quelques points d'alg{\`e}bre homologique
et d'alg{\`e}bre homotopique et leurs applications}, Proceedings
of Tbilisi A.\ Razmadze Mathematical Institute \textbf{48} (1975),
5--137, in Russian.

\bibitem{Inasaridze70}
H.\ Inassaridze, \emph{Cohomology with values in semigroups},
Ph.D.\ thesis abstract, Tbilisi, 1970, in Russian.

\bibitem{Dold-Puppe}
A.~Dold and D.~Puppe, \emph{Homologie nicht-additiver {F}unktoren.
  {A}nwendungen}, Ann.\ Inst.\ Fourier (Grenoble) \textbf{11} (1961), 201--312.

\bibitem{Eilenberg-Moore}
S.~Eilenberg and J.~Moore, \emph{Foundations of relative homological algebra}, Memoirs AMS, vol.~55, American Mathematical Society, 1965.

\bibitem{Tierney-Vogel2}
M.~Tierney and W.~Vogel, \emph{Simplicial resolutions and derived functors}, Math.\ Z.\ \textbf{111} (1969), no.~1, 1--14.

\bibitem{Swan}
R.~G. Swan, \emph{Some relations between higher {K}-functors}, J.~Algebra \textbf{21} (1972), 113--136.

\end{thebibliography}
\end{document}